# A concise modification of Marshall-Olkin family of distributions for reliability analysis

**Mueen-ud-Din Azad[1,*] and Muhammad Mohsin[2]**

[1] Department of Quantitative Methods, University of Management and Technology, Lahore-PAKISTAN
[2] Department of Statistics, COMSATS University Islamabad, Lahore Campus, Lahore-PAKISTAN

*Corresponding author, email: mueen.azad@umt.edu.pk

**Abstract:**

The significance of Marshall-Olkin distribution in reliability theory has motivated us to introduce a generalized exponentiated Marshall-Olkin (GEMO) family of distributions along with its five sub models. These sub models are GEMO–Exponential, GEMO–Weibull, GEMO–Gamma, GEMO–Lomax and GEMO–Log Normal distributions. The proposed GEMO distribution has high probability at right tail which is quite helpful to deal with the extreme events in reliability data. Several reliability characteristics of the proposed distribution are studied. Two real life examples are presented to demonstrate the usefulness of GEMO-Weibull distribution. The comparison of GEMO-Weibull distribution with the existing life time distributions is also given which endorses the performance of our model.

## 1. Introduction

Reliability theory, a multidisciplinary science, has a significant role to measure the consistency of a system and has close links with the probability distributions. Advancement in science and technology makes us come upon more powerful, complex and sophisticated engineering systems. These high tech engineering systems are likely to face more sophisticated reliability issues which ultimately require more complicated models for their solutions (Kuo and Zuo 2003). Probability distributions provide such models that are commonly and effectively used to measure and analyze the reliability of a system.

Many reliability models such as life time, survival time, k-out-of-n, stress and strength, multistate and maintainable models have been presented in literature. A brief summary of the most significant work includes: applications to the yield strength and fatigue life of steel Weibull (1951), fracture strength of glass Keshavan (1980), pitting corrosion in pipes Sheikh et al. (1990), adhesive wear in metals Qureshi and Sheikh (1997), failure of carbon fiber composites Durham and Padgett (1997), coatings Almeida, (1999), brittle materials Fok, et al. (2001), composite materials Newell, et al. (2002), concrete components Li, et al. (2003), stochastic aging and dependence for reliability Lai and Xie (2006), game theory to explore the nature of optimal investments in the security of simple series and parallel systems Azaiez and Bier (2007), Systems requiring very high levels of reliability, such as aircraft controls or spacecraft Myers (2007) and a new parametric distribution for modeling lifetime data with bathtub shaped hazard rate Lai et al. (2016).

In recent years, researchers use different base line distribution to develop Marshall- Olkin (1997) family of distributions. For reference, Saboor and Pogány (2016) introduce a Marshall–Olkin variant of the Provost Type gamma–Weibull probability distribution. Nojosa and Rathie (2017) develop the Marshall-Olkin Rathie Swamee distribution and study its distributional properties. Rocha et al. (2017) propose a new way to generate defective distributions to model cure fractions derived from the Marshall Olkin family of distributions. Afify et al. (2018) study the Marshall-Olkin additive Weibull distribution in order to allow a wide variation in the shape of the hazard rate, including increasing, decreasing, bathtub and unimodal shapes. Ngelova et al. (2019) study the characteristic of the cummulative function of Marshall-Olkin Inverse Lomax distribution to the horizontal asymptote with respect to the Hausdorff distance. Muhammad and Liu (2019) make the characterization of the Marshall-Olkin-G family of distribution by left and right truncated moments. Conleth et al. (2019) propose a new Marshall-Olkin extended Weibull-exponential distribution. Okasha et al. (2020) introduce a new model of lifetime distributions which is called Marshall-Olkin extended inverse Weibull distribution. Cui et al. (2020) propose a new Marshall-Olkin Weibull Distribution and discuss its importance in reliability.

Marshall and Olkin (1997) introduce a class of distributions to add a parameter into a family of baseline distributions. They define one-parameter family of survival functions as

$$\bar{G}(x;\alpha) = \frac{\alpha \bar{F}(x)}{1 - \bar{\alpha}\bar{F}(x)},$$

where $\bar{F}(x)$ is a survival function of the baseline distribution and $\alpha$ is an additional parameter. They also propose two-parameters exponential and three-parameters Weibull distributions in their pioneer work. For a comprehensive discussion on life distributions, see Marshall and Olkin (2007). Jose (2011) shows the application of Marshall-Olkin technique on various areas like reliability theory, time series modeling and stress and strength analysis. Cordeiro et al. (2014) propose three new distributions on the basis of the Marshall-Olkin scheme and discuss their mathematical and statistical properties as well. Dias et al. (2016) develop a new class with three extra shape parameters and call it Exponentiated Marshall-Olkin (EMO-G) family of distributions. Handique et al. (2019) propose the exponentiated generalized Marshall–Olkin (EGMO) family of distributions with three additional parameters and also prove the mathematical and statistical properties of the proposed family.

In various applications of reliability we face difficulties to model skewed and high or low kurtosis data. Although some lifetime distributions such as exponential, Weibull and gamma distributions have been extensively used in reliability analysis but they have certain limitations like, exponential distribution has constant hazard rate. Due to such limitations in the baseline distributions, the researchers have always been generating new distributions or modifying the existing distributions to manipulate such difficulties. In this paper we made effort to modify the Marshall-Olkin family of distributions by introducing a new generator which shows versatile aptitude to model the real situations more effectively. The proposed family of distributions having two parameters successfully models skewed and high or low kurtosis data. The comparison of some important baseline distributions such as Exponentiated Marshall-Olkin (EMO-G) family of distributions, Weibull, exponential and gamma with our proposed family of distributions clearly asserts its compatibility for wide range of applications.

The article is organized as follows: In Section 2, we use a new generator in Marshall-Olkin distribution to develop another family of distributions and study some of its distributional properties. The statistical properties regarding reliability analysis are discussed in Section 3. Estimation of the model parameters, an asymptotic confidence interval, and likelihood ratio test are discussed in Section 4. To establish the adequacy of the model, two real life examples regarding reliability are presented along with its comparison with some baseline distributions in Section 5. Finally, concluding remarks are presented in Section 6.

## 2. The Proposed Family and its Characteristics

In this section, we propose a new family of lifetime distributions named as generalized exponentiated Marshall-Olkin (GEMO) family of distributions. The cummulative distribution function (cdf), probability density function (pdf), survival function (sf) and hazard rate function (hrf) are given by

$$G(x;\alpha,\beta,\gamma,\xi) = 1 - \left[\frac{\alpha \bar{F}^{\gamma}(x)}{1+(\bar{\alpha})\bar{F}^{\gamma}(x)}\right]^{\beta}, \qquad (1)$$

$$g(x;\alpha,\beta,\gamma,\xi) = \frac{\beta\gamma\alpha^{\beta} f(x)\left[1-(1-F(x))^{\gamma}(\bar{\alpha})\right]^{-(\beta+1)}\left[1-F(x)\right]^{\gamma\beta}}{1-F(x)}, \quad I(-\infty < x < \infty), \qquad (2)$$

$$\bar{G}(x;\alpha,\beta,\gamma,\xi) = \left[\frac{\alpha \bar{F}^{\gamma}(x)}{1+(\bar{\alpha})\bar{F}^{\gamma}(x)}\right]^{\beta}, \qquad (3)$$

$$h(x;\alpha,\beta,\gamma,\xi) = \frac{\beta\gamma\alpha^{\beta} f(x)\left[1-(1-F(x))^{\gamma}(\bar{\alpha})\right]^{-(\beta+1)}\left[1-F(x)\right]^{\gamma\beta}}{\left[1-F(x)\right]\left[\left[\frac{\alpha \bar{F}^{\gamma}(x)}{1+(\bar{\alpha})\bar{F}^{\gamma}(x)}\right]^{\beta}\right]}, \qquad (4)$$

respectively where $\bar{F}(x)$ and $f(x)$ are the sf and pdf of the baseline distribution and $\alpha > 0, \beta > 0, \gamma > 0$ are the additional scale and the shape parameters. The $\xi$ is the parameter vector of the baseline distribution. In addition $G(x) = F(x)$, if $\alpha = \beta = \gamma = 1$.

By using the generalized binomial expression, we can show that the pdf (2) of $X$ confesses the following expression

$$g(x) = \beta\gamma\alpha^{\beta} f(x) \sum_{j=0}^{\infty} w_j \left[1-F(x)\right]^{\gamma(\beta+j)-1}, \qquad (5)$$

where $w_j = w_j^* (-1)^j (1-\alpha)^j$, $w_j^* = \dfrac{(\beta+j)!}{\beta!\, j!}$.

The new proposed family of distributions is applicable to all the baseline functions that satisfy the properties of $F(x)$. Below are some special distributions derived from the proposed family of distributions.

### a. The GEMO–Exponential (GEMO–E) Distribution

Suppose the exponential distribution with parameter $\lambda > 0$, $f(x) = \lambda e^{-\lambda x}$ and $F(x) = 1 - e^{-\lambda x}$ be the pdf and cdf of the baseline distribution, then for $GEMO-E(\alpha, \beta, \gamma, \lambda)$ the cdf and pdf will be

$$G_{GEMO-E}(x; \alpha, \beta, \gamma, \lambda) = 1 - \left( \frac{\alpha e^{-x\lambda}}{1 + (1-\alpha)\left(e^{-x\lambda}\right)^\gamma} \right)^\beta \text{ and}$$

$$g_{GEMO-E}(x; \alpha, \beta, \gamma, \lambda) = \alpha^\beta \beta \gamma \lambda \left(e^{-x\lambda}\right)^{\beta\gamma} \left(1 - (1-\alpha)\left(e^{-x\lambda}\right)^\gamma\right)^{-(\beta+1)}$$

respectively. Figure 1 (a-c) shows the graphical presentation of the pdf for GEMO-Exponential distribution. In these figures we note the role of the three additional parameters.

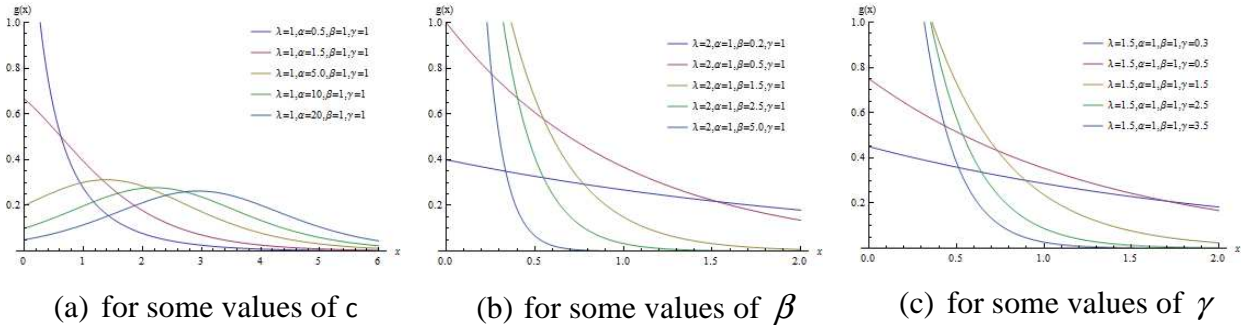

(a) for some values of c  (b) for some values of $\beta$  (c) for some values of $\gamma$

**Figure 1:** pdf plots for GEMO-Exponential distribution for different values of $\alpha, \beta$ and $\gamma$.

### b. The GEMO–Weibull (GEMO–W) Distribution

Suppose the Weibull distribution with parameter $\lambda, \theta > 0$, $f(x) = (\lambda/\theta)(x/\theta)^{\lambda-1} e^{-(x/\theta)^\lambda}$ and $F(x) = 1 - e^{-(x/\theta)^\lambda}$ be the pdf and cdf of the base line distribution, then for $GEMO-W(\alpha, \beta, \gamma, \lambda, \theta)$ the cdf and pdf will be

$$G_{GEMO-W}(x;\alpha,\beta,\gamma,\lambda,\theta)=1-\left[\alpha\left(e^{-\left(\frac{x}{\theta}\right)^{\lambda}}\right)^{\gamma}\left(1+(1-\alpha)\left(e^{-\left(\frac{x}{\theta}\right)^{\lambda}}\right)^{\gamma}\right)^{-1}\right]^{\beta}$$ and

$$g_{GEMO-W}(x;\alpha,\beta,\gamma,\lambda,\theta)=\frac{\lambda\beta\gamma\alpha^{\beta}}{\theta}\left(\frac{x}{\theta}\right)^{\lambda-1}\left(e^{-\left(\frac{x}{\theta}\right)^{\lambda}}\right)^{\beta\gamma}\left(1-\left(e^{-\left(\frac{x}{\theta}\right)^{\lambda}}\right)^{\gamma}(1-\alpha)\right)^{-(\beta+1)}$$

respectively. Figure 2 (a-c) is the graphical presentation of the pdf for GEMO-Weibull distribution. In these figures we note the role of the three additional parameters. These additional parameters help us to model skewed and high or low kurtosis data.

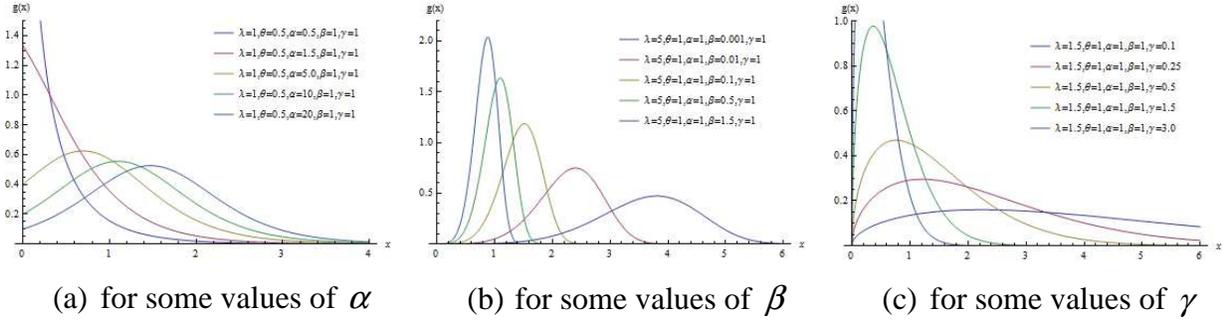

(a) for some values of $\alpha$     (b) for some values of $\beta$     (c) for some values of $\gamma$

**Figure 2:** pdf plots for GEMO-Weibull distribution for different values of $\alpha, \beta$ and $\gamma$.

### c. The GEMO–Gamma (GEMO–G) Distribution

Suppose the Gamma distribution with parameter $\lambda, \theta > 0$, $f(x)=\left(\theta^{\lambda}/\Gamma(\lambda)\right)x^{\lambda-1}e^{-\theta x}$ and $F(x)=(1/\Gamma(\lambda))(\Gamma_{x}(\lambda,\theta x))$ be the pdf and cdf of the base line distribution, then for $GEMO-G(\alpha,\beta,\gamma,\lambda,\theta)$ the cdf and pdf will be

$$G_{GEMO-G}(x;\alpha,\beta,\gamma,\lambda,\theta)=1-\left(\frac{\alpha\left(1-\frac{\Gamma_{x}(\lambda,x\theta)}{\Gamma(\lambda)}\right)}{1+(1-\alpha)\left(1-\frac{\Gamma_{x}(\lambda,x\theta)}{\Gamma(\lambda)}\right)^{\gamma}}\right)^{\beta}$$ and

$$g_{GEMO-G}(x;\alpha,\beta,\gamma,\lambda,\theta)=\frac{\alpha^{\beta}\gamma\beta\theta^{\lambda}e^{-x\theta}x^{\lambda-1}}{\Gamma(\lambda)}\left(1-\frac{\Gamma_{x}(\lambda,x\theta)}{\Gamma(\lambda)}\right)^{\beta\gamma-1}\left(1-(1-\alpha)\left(1-\frac{\Gamma_{x}(\lambda,x\theta)}{\Gamma(\lambda)}\right)^{\gamma}\right)^{-(\beta+1)}$$

respectively. Figure 3 (a-c) depicts the graphical presentation of the pdf for GEMO-Gamma distribution. In these figures we note the role of the three additional parameters. These additional parameters help us to model skewed and heavy tailed data.

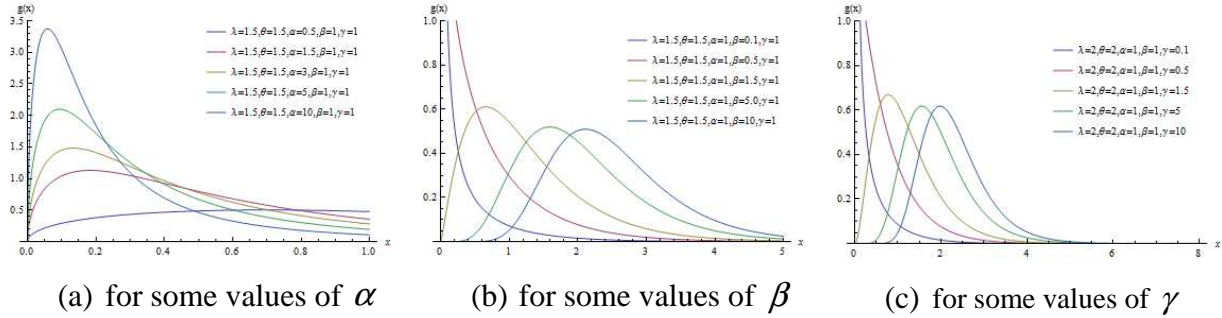

(a) for some values of $\alpha$      (b) for some values of $\beta$      (c) for some values of $\gamma$

**Figure 3:** pdf plots for GEMO-Gamma distribution for different values of $\alpha, \beta$ and $\gamma$.

### d. The GEMO–Lomax (GEMO–L) Distribution

Suppose the Lomax distribution with parameter $\lambda, \theta > 0$, $f(x) = (\theta/\lambda)(1 + x/\lambda)^{-(\theta+1)}$ and $F(x) = 1 - (1 + x/\lambda)^{-\theta}$ be the pdf and cdf of the base line distribution, then for $GEMO - L(\alpha, \beta, \gamma, \lambda, \theta)$ the cdf and pdf will be

$$G_{GEMO-L}(x; \alpha, \beta, \gamma, \lambda, \theta) = 1 - \left( \frac{\alpha \left(1 + \frac{x}{\lambda}\right)^{-\theta}}{1 + (1-\alpha)\left(1 + \frac{x}{\lambda}\right)^{-\theta\gamma}} \right)^{\beta}$$

and

$$g_{GEMO-L}(x; \alpha, \beta, \gamma, \lambda, \theta) = \frac{\theta \beta \gamma \alpha^{\beta}}{\lambda} \frac{\left(1 + \frac{x}{\lambda}\right)^{-(\theta+1)} \left(\left(1 + \frac{x}{\lambda}\right)^{-\theta}\right)^{\beta\gamma - 1}}{\left(1 - (1-\alpha)\left(1 + \frac{x}{\lambda}\right)^{-\theta\gamma}\right)^{\beta+1}}$$

respectively. Figure 4 (a-c) shows the graphical presentation of the pdf for GEMO-Lomax distribution. In these figures we note the role of the three additional parameters.

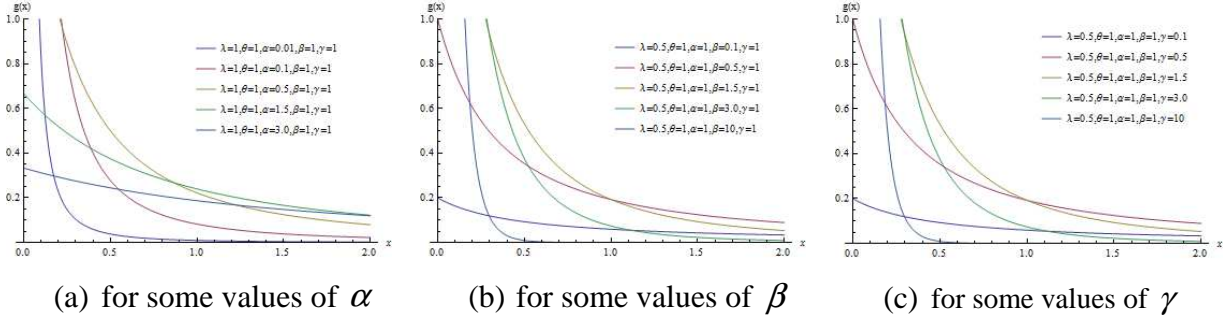

(a) for some values of $\alpha$     (b) for some values of $\beta$     (c) for some values of $\gamma$

**Figure 4:** pdf plots for GEMO-Lomax distribution for different values of $\alpha, \beta$ and $\gamma$.

### e. The GEMO–Log Normal (GEMO–LN) Distribution

Suppose the Log-Normal distribution with parameter $\mu, \sigma > 0$, $f(x) = (1/\sigma x\sqrt{2\pi})e^{-\left((\ln(x)-\mu)^2/2\sigma^2\right)}$ and $F(x) = (1/2) + (1/2)\,Erf\left((-\mu + \ln(x))/\sqrt{2}\sigma\right)$ be the pdf and cdf of the base line distribution, then for $GEMO-LN(\alpha, \beta, \gamma, \mu, \sigma)$ the cdf and pdf will be

$$G_{GEMO-LN}(x;\alpha,\beta,\gamma,\mu,\sigma) = 1 - \left(\frac{\alpha\left(\frac{1}{2} - \frac{1}{2}erf\left(\frac{-\mu + Log(x)}{\sqrt{2}\sigma}\right)\right)}{1 + (1-\alpha)\left(\frac{1}{2} - \frac{1}{2}erf\left(\frac{-\mu + Log(x)}{\sqrt{2}\sigma}\right)\right)^\gamma}\right)^\beta$$

and

$$g_{GEMO-LN}(x;\alpha,\beta,\gamma,\mu,\sigma) = \frac{\gamma\beta\alpha^\beta e^{-\frac{(-\mu+Log(x))^2}{2\sigma^2}}}{\sqrt{2\pi}\,x\sigma} \cdot \frac{\left(\left(\frac{1}{2} - \frac{1}{2}erf\left(\frac{-\mu + Log(x)}{\sqrt{2}\sigma}\right)\right)\right)^{\beta\gamma-1}}{\left(1 - (1-\alpha)\left(\frac{1}{2} - \frac{1}{2}erf\left(\frac{-\mu + Log(x)}{\sqrt{2}\sigma}\right)\right)^\gamma\right)^{\beta+1}}$$

respectively. Figure 5 (a-c) provides the graphical presentation of the pdf for GEMO-Log Normal distribution. In these figures we note the role of the three additional parameters.

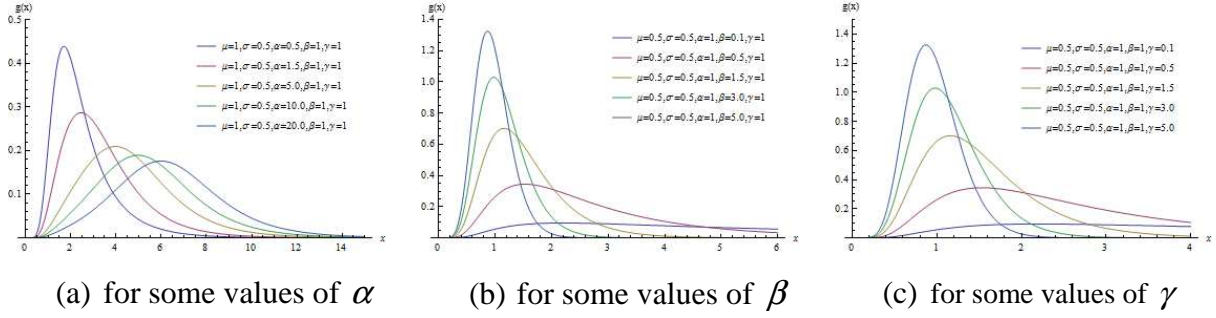

(a) for some values of $\alpha$  (b) for some values of $\beta$  (c) for some values of $\gamma$

**Figure 5:** pdf plots for GEMO-Log Normal distribution for different values of $\alpha, \beta$ and $\gamma$.

The probability weighted moment (PWM), introduced by Greenwood et al. (1979), is the expected value of a function of a random variable whose mean exists. The $(l, j, k)$ th PWM of $X$ having cdf $F(x)$ is defined as

$$M_{l,j,k} = \int_{-\infty}^{\infty} x^l \left[F(x)\right]^j \left[1-F(x)\right]^k f(x)\, dx.$$

**Theorem 2.1**

If a random variable $X$ follows the GEMO, then the moment generating function $M(t)$ is given by

$$M(t) = \gamma\beta\alpha^\beta \sum_{j=0}^{\infty} w_j \sum_{r=0}^{\infty} \frac{t^r}{r!} E\left[x^r \left(\bar{F}(x)\right)^{\gamma(\beta+j)-1}\right]. \tag{6}$$

**Corollary 2.1**

If a random variable $X$ follows the GEMO, then the r$^{th}$ moment $E(X^r)$ is

$$\mu'_r = \beta\gamma\alpha^\beta \sum_{j=0}^{\infty} w_j E\left[x^r \left(1-F(x)\right)^{\gamma(\beta+j)-1}\right]. \tag{7}$$

**Theorem 2.2**

If a random variable $X$ follows the GEMO, then the quantile function $Q = G^{-1}(x; \alpha, \beta, \gamma, \xi)$ is

$$x = G^{-1}\left(1 - \left\{\frac{(1-u)^{1/\beta}}{\alpha + \bar{\alpha}(1-u)^{1/\beta}}\right\}^{1/\gamma}\right). \tag{8}$$

## 2.1 Entropy

Entropy of a random variable $X$ is the measure of variation of uncertainty in a given population. In the literature different entropy measures are available. In this sub-section we discuss Varma (1966) entropy and Shannon (1948) entropy.

**Theorem 2.3**

If a random variable $X$ follows the GEMO, then the Varma entropy $H_v(a,b)$ is given by

$$H_v(a,b) = \frac{1}{b-a} \log \int [g(x)]^{a+b-1} dx, \qquad \text{for } b-1 < a < b; \ b \geq 1, \qquad (9)$$

where $g(x)$ is defined in (2). Shannon entropy can be derived from the above equation when $a \to 1, b \to 1$.

## 2.2 Distribution of Order Statistics

Order statistics has a significant role in robust location estimation, reliability, detection of outliers, censored sampling, strength of materials, quality control, selecting the best and many more, Arnold et al, (2008).

**Theorem 2.4**

If a random variable $X$ follows the GEMO, then the r$^{th}$ order statistics $g_{r:n}(x)$ is given by

$$f_{r:n}(x) = \frac{g(x)}{B(r, n-r+1)} G^{r-1}(x) [1 - G(x)]^{n-r}, \qquad (10)$$

where $g(x)$ and $G(x)$ are the pdf and cdf of GEMO family of distribution respectively.

## 3. Reliability Analysis

In this section, we discuss some reliability properties like mean residual life (MRL) function, mean past lifetime (MPL) and conditional moments of GEMO.

### 3.1 Survival and Hazard Functions

The survival function (sf) and the hazard rate function (hrf) for GEMO family of distributions are defined in equation (3) and (4) respectively.

The survival and hazard rate functions for GEMO-Exponential are expressed as

$$S_{GEMO-E}(t;\alpha,\beta,\gamma,\lambda) = \left(\frac{\alpha e^{-t\lambda}}{1+(1-\alpha)\left(e^{-t\lambda}\right)^{\gamma}}\right)^{\beta} \text{ and}$$

$$h_{GEMO-E}(t;\alpha,\beta,\gamma,\lambda) = \frac{\alpha^{\beta}\beta\gamma\lambda\left(e^{-t\lambda}\right)^{\beta\gamma}}{\left(1-(1-\alpha)\left(e^{-t\lambda}\right)^{\gamma}\right)^{-(\beta+1)}} \left(\frac{\alpha e^{-t\lambda}}{1+(1-\alpha)\left(e^{-t\lambda}\right)^{\gamma}}\right)^{-\beta}$$

respectively. Figure 6 (a-c) shows the graphical presentation of the hrf for GEMO-Exponential distribution. In these figures we can note the role of the three additional parameters. These additional parameters help us to model increasing, decreasing and, bath tub shaped hazard rate function.

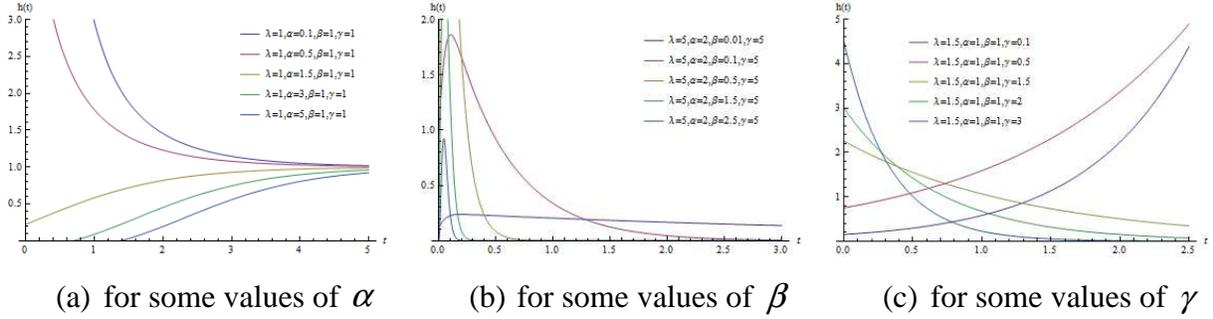

(a) for some values of $\alpha$      (b) for some values of $\beta$      (c) for some values of $\gamma$

**Figure 6:** hrf plots for GEMO-Exponetial distribution for different values of $\alpha, \beta$ and $\gamma$.

The survival and hazard rate functions for GEMO-Weibull are given as

$$S_{GEMO-W}(t;\alpha,\beta,\gamma,\lambda,\theta) = \left[\alpha\left(e^{-\left(\frac{t}{\theta}\right)^{\lambda}}\right)^{\gamma}\left(1+(1-\alpha)\left(e^{-\left(\frac{t}{\theta}\right)^{\lambda}}\right)^{\gamma}\right)^{-1}\right]^{\beta} \text{ and}$$

$$h_{GEMO-W}(t;\alpha,\beta,\gamma,\lambda,\theta) = \frac{\lambda\beta\gamma\alpha^{\beta}}{\theta} \frac{\left(\frac{t}{\theta}\right)^{\lambda-1}\left(e^{-\left(\frac{t}{\theta}\right)^{\lambda}}\right)^{\beta\gamma}}{\left(1-\left(e^{-\left(\frac{t}{\theta}\right)^{\lambda}}\right)^{\gamma}(1-\alpha)\right)^{(\beta+1)}} \left(\frac{\alpha\left(e^{-\left(\frac{t}{\theta}\right)^{\lambda}}\right)^{\gamma}}{1+(1-\alpha)\left(e^{-\left(\frac{t}{\theta}\right)^{\lambda}}\right)^{\gamma}}\right)^{-\beta}$$

respectively. Figure 7 (a-c) shows the graphical presentation of the hrf for GEMO-Weibull distribution. In these figures we can note the role of the three additional parameters. These additional parameters help us to model increasing, decreasing and, bath tub shaped hazard rate function.

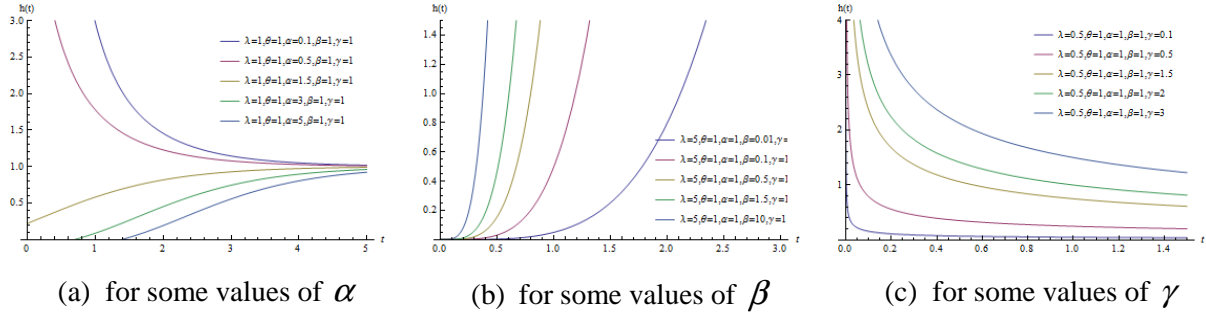

(a) for some values of $\alpha$   (b) for some values of $\beta$   (c) for some values of $\gamma$

**Figure 7:** hrf plots for GEMO-Weibull distribution for different values of $\alpha, \beta$ and $\gamma$.

The survival and hazard rate functions for GEMO-Gamma are

$$S_{GEMO-G}(t;\alpha,\beta,\gamma,\lambda,\theta) = \left( \frac{\alpha\left(1-\frac{\Gamma_x(\lambda,t\theta)}{\Gamma(\lambda)}\right)}{1+(1-\alpha)\left(1-\frac{\Gamma_x(\lambda,t\theta)}{\Gamma(\lambda)}\right)^\gamma} \right)^\beta \quad \text{and}$$

$$h_{GEMO-G}(t;\alpha,\beta,\gamma,\lambda,\theta) = \frac{\alpha^\beta \gamma\beta\theta^\lambda e^{-x\theta} t^{\lambda-1}}{\Gamma(\lambda)} \left(1-\frac{\Gamma_t(\lambda,t\theta)}{\Gamma(\lambda)}\right)^{\beta\gamma-1} \left(1-(1-\alpha)\left(1-\frac{\Gamma_t(\lambda,t\theta)}{\Gamma(\lambda)}\right)^\gamma\right)^{-(\beta+1)}$$

$$\left( \frac{\alpha\left(1-\frac{\Gamma_t(\lambda,t\theta)}{\Gamma(\lambda)}\right)}{1+(1-\alpha)\left(1-\frac{\Gamma_t(\lambda,t\theta)}{\Gamma(\lambda)}\right)^\gamma} \right)^{-\beta}$$

respectively. Figure 8 (a-c) gives the graphical presentation of the hrf for GEMO-Gamma distribution. In these figures we can note the role of the three additional parameters. These additional parameters help us to model increasing, decreasing and, bath tub shaped hazard rate functions.

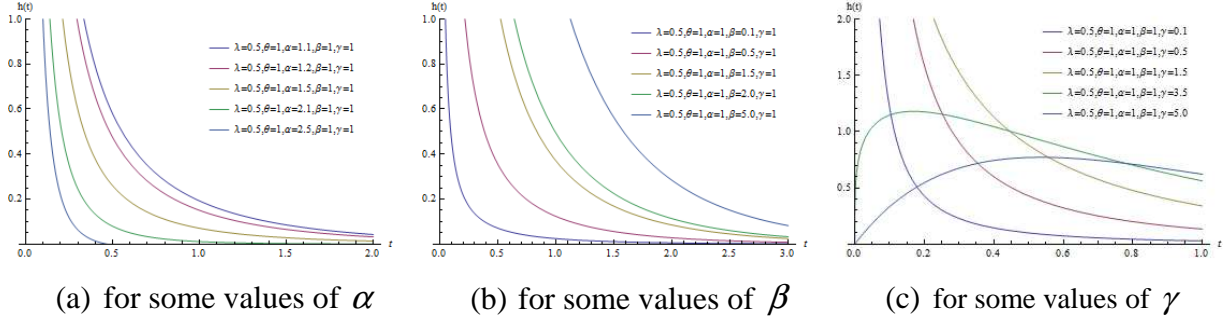

(a) for some values of $\alpha$  (b) for some values of $\beta$  (c) for some values of $\gamma$

**Figure 8:** hrf plots for GEMO-Gamma distribution for different values of $\alpha, \beta$ and $\gamma$.

The survival and hazard rate functions for GEMO-Lomax are given as

$$S_{GEMO-L}(t;\alpha,\beta,\gamma,\lambda,\theta) = \left(\frac{\alpha\left(1+\frac{t}{\lambda}\right)^{-\theta}}{1+(1-\alpha)\left(1+\frac{t}{\lambda}\right)^{-\theta\gamma}}\right)^{\beta} \quad \text{and}$$

$$h_{GEMO-L}(t;\alpha,\beta,\gamma,\lambda,\theta) = \frac{\theta\beta\gamma\alpha^{\beta}}{\lambda} \frac{\left(1+\frac{t}{\lambda}\right)^{-(\theta+1)}\left(\left(1+\frac{t}{\lambda}\right)^{-\theta}\right)^{\beta\gamma-1}}{\left(1-(1-\alpha)\left(1+\frac{t}{\lambda}\right)^{-\theta\gamma}\right)^{\beta+1}}\left(\frac{\alpha\left(1+\frac{t}{\lambda}\right)^{-\theta}}{1+(1-\alpha)\left(1+\frac{t}{\lambda}\right)^{-\theta\gamma}}\right)^{-\beta}$$

respectively. Figure 9 (a-c) shows the graphical presentation of the hrf for GEMO-Lomax distribution. In these figures we can note the role of the three additional parameters. These additional parameters help us to model increasing, decreasing and, bath tub shaped hazard rate function

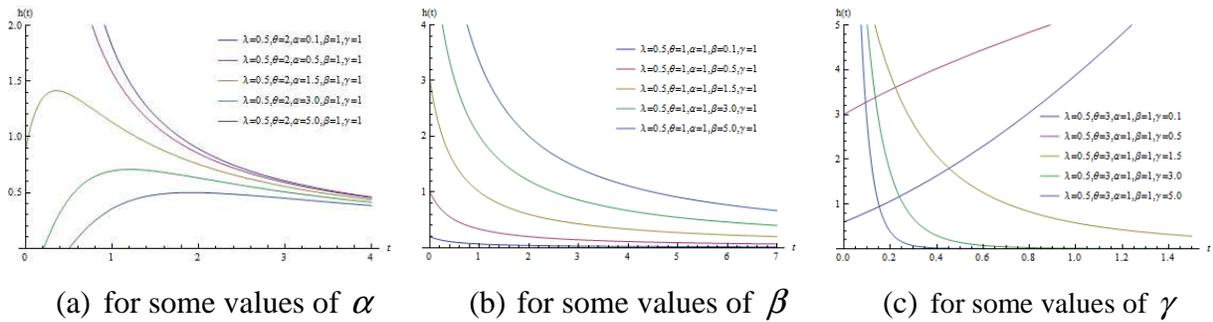

(a) for some values of $\alpha$  (b) for some values of $\beta$  (c) for some values of $\gamma$

**Figure 9:** hrf plots for GEMO-Lomax distribution for different values of $\alpha, \beta$ and $\gamma$.

The survival and hazard rate functions for GEMO-Log Normal are

$$S_{GEMO-LN}(t;\alpha,\beta,\gamma,\mu,\sigma) = 1 - \left( \frac{\alpha\left(\frac{1}{2} - \frac{1}{2}erf\left(\frac{-\mu + Log(t)}{\sqrt{2}\sigma}\right)\right)}{1 + (1-\alpha)\left(\frac{1}{2} - \frac{1}{2}erf\left(\frac{-\mu + Log(t)}{\sqrt{2}\sigma}\right)\right)^{\gamma}} \right)^{\beta} \quad \text{and}$$

$$h_{GEMO-LN}(t;\alpha,\beta,\gamma,\mu,\sigma) = \frac{\gamma\beta\alpha^{\beta}e^{-\frac{(-\mu+Log(t))^2}{2\sigma^2}}}{\sqrt{2\pi}\,x\sigma} \frac{\left(\left(\frac{1}{2} - \frac{1}{2}erf\left(\frac{-\mu + Log(t)}{\sqrt{2}\sigma}\right)\right)\right)^{\beta\gamma-1}}{\left(1 - (1-\alpha)\left(\frac{1}{2} - \frac{1}{2}erf\left(\frac{-\mu + Log(t)}{\sqrt{2}\sigma}\right)\right)^{\gamma}\right)^{\beta+1}}$$

$$\left( \frac{\alpha\left(\frac{1}{2} - \frac{1}{2}erf\left(\frac{-\mu + Log(t)}{\sqrt{2}\sigma}\right)\right)}{1 + (1-\alpha)\left(\frac{1}{2} - \frac{1}{2}erf\left(\frac{-\mu + Log(t)}{\sqrt{2}\sigma}\right)\right)^{\gamma}} \right)^{-\beta}$$

respectively. Figure 10 (a-c) shows the graphical presentation of the hrf for GEMO-Log Normal distribution. In these figures we can note the role of the three additional parameters. These additional parameters help us to model increasing, decreasing and, bath tub shaped hazard rate function

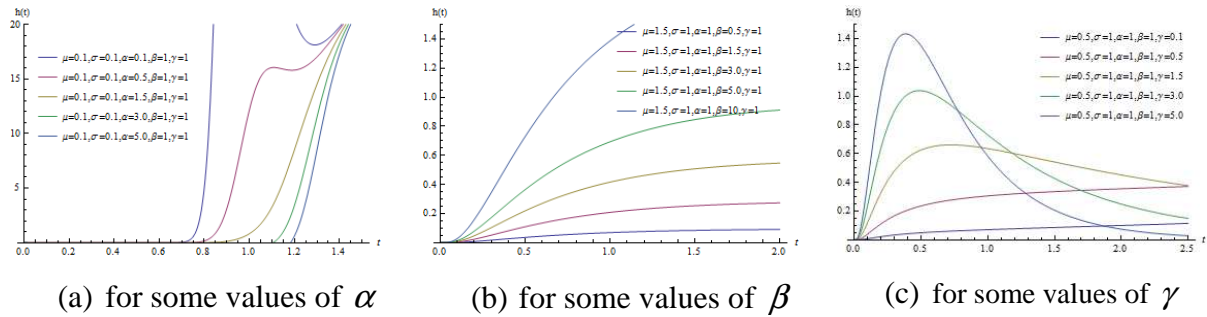

(a) for some values of $\alpha$    (b) for some values of $\beta$    (c) for some values of $\gamma$

**Figure 10:** hrf plots for GEMO-Log Normal distribution for different values of $\alpha, \beta$ and $\gamma$.

The mean residual life (MRL) function plays an important role, as it measures the expected value of the remaining lifetime of the system provided that it has survived up to time $t$.

**Theorem 3.1**

If a random variable $X$ follows the GEMO then the mean residual life function $\mu(t)$ is given by

$$\mu(t) = \beta\gamma\alpha^{-\beta}\sum_{j=0}^{\infty} w_j \int_t^{\infty} x\ f(x)[1-F(x)]^{\gamma(\beta+j)-1}\ dx - t. \qquad (11)$$

**Proof:** The MRL function is the expected remaining life of $(T-t)$, given that the item has survived to time $t$ and defined as

$$\mu(t) = E[(X-t)\,|\,X > t] = \frac{\int_t^{\infty} x\, g(x)\, dx}{1 - G(t)} - t,$$

$$\mu(t) = \left\{ \frac{\beta\gamma\alpha^{\beta}\beta \int_t^{\infty} x\ f(x) \sum_{j=0}^{\infty} w_j [1-F(x)]^{\gamma(\beta+j)-1}\ dx}{\left[\alpha(1-F(t))^{\gamma}\right]^{\beta}\left[1-(1-\alpha)(1-F(t))^{\gamma}\right]^{-\beta}} \right\} - t,$$

$$\mu(t) = \beta\gamma\alpha^{-\beta}\sum_{j=0}^{\infty} w_j \int_t^{\infty} x\ f(x)[1-F(x)]^{\gamma(\beta+j)-1}\ dx - t.$$

The mean past life time (MPL) corresponds to the mean time elapsed since the failure of $X$ given that $X \leq t$. It is used to access the lifetime of the system.

**Theorem 3.2**

If a random variable $T$ follows GEMO, then the mean past lifetime $k(t)$ is given by,

$$k(t) = t - \frac{\beta\gamma\alpha^{\beta}\left[1+(\alpha-1)(1-F(t))^{\gamma}\right]^{\beta}}{\left[(1-F(t))^{\gamma}\right]^{\beta}(1+(\alpha-1)-\alpha)^{\beta}} \sum_{j=0}^{\infty} w_j \int_0^t x f(x)[1-F(x)]^{\gamma(\beta+j)-1}\ dx. \qquad (12)$$

**Proof:** The MPL which measures the expected value of the remaining lifetime of the system provided that it has survived up to time $t$ is given by

$$k(t) = E[t - X\,|\,X \leq t] = \frac{\int_0^t G(x)dx}{G(t)} = t - \frac{\int_0^t x\, g(x)\, dx}{G(t)},$$

$$k(t)=t-\frac{\beta\gamma\alpha^{\beta}\left[1+(\alpha-1)(1-F(t))^{\gamma}\right]^{\beta}}{\left[(1-F(t))^{\gamma}\right]^{\beta}(1+(\alpha-1)-\alpha)^{\beta}}\sum_{j=0}^{\infty}w_{j}\int_{0}^{t}xf(x)[1-F(x)]^{\gamma(\beta+j)-1}\,dx.$$

**Theorem 3.3**

If a random variable $X$ follows GEMO then conditional moment $E(X^n \mid X \geq t)$ is given by

$$E(X^n \mid X \geq t) = \left[\frac{\alpha(1-F(t))^{\gamma}}{1-(1-\alpha)(1-F(t))^{\gamma}}\right]^{-\beta}\sum_{j=0}^{\infty}w_{j}\int_{t}^{\infty}x^{n}f(x)[1-F(x)]^{\gamma(\beta+j)-1}\,dx. \tag{13}$$

**Proof:** By definition the conditional moments are

$$E(X^n \mid X \geq t) = \frac{1}{S(t)}\int_{t}^{\infty}x^{n}f(x)\,dx,$$

$$E(X^n \mid X \geq t) = \left[\frac{\alpha(1-F(t))^{\gamma}}{1-(1-\alpha)(1-F(t))^{\gamma}}\right]^{-\beta}\sum_{j=0}^{\infty}w_{j}\int_{t}^{\infty}x^{n}f(x)[1-F(x)]^{\gamma(\beta+j)-1}\,dx.$$

## 4. Estimation

The likelihood function for the distribution proposed in (2) on the basis of $n$ independent observations takes the form

$$L(\underline{x}\mid\alpha,\beta,\gamma)=\alpha^{n\beta}\beta^{n}\gamma^{n}\prod_{i=1}^{n}[f(x_{i})]\prod_{i=1}^{n}[1-F(x_{i})]^{\beta\gamma-1}\prod_{i=1}^{n}\left[1+(\alpha-1)(1-F(x_{i}))^{\gamma}\right]^{-(\beta+1)}.$$

The log-likelihood function is presented by

$$l = \ln L(\underline{x}\mid\alpha,\beta,\gamma) = n\ln(\beta) + n\beta\ln(\alpha) + n\ln(\gamma) + \sum_{i=1}^{n}\ln(f(x_{i})) + (\beta\gamma-1)\sum_{i=1}^{n}\ln(1-F(x_{i}))$$

$$-(\beta+1)\sum_{i=1}^{n}\ln\left(1+(\alpha-1)(1-F(x_{i}))^{\gamma}\right).$$

To find the maximum likelihood estimates of $(\hat{\alpha},\hat{\beta},\hat{\gamma})^{T}$ we compute the corresponding estimating equations by taking partial derivative of the log-likelihood function with respect to $(\alpha,\beta,\gamma)^{T}$ and equating them to zero:

$$\frac{\partial[\ln L(\underline{x}\mid\alpha,\beta,\gamma)]}{\partial\alpha} = \frac{n\beta}{\alpha} - (\beta+1)\sum_{i=1}^{n}\frac{(1-F(x_{i}))^{\gamma}}{1-(1-\alpha)(1-F(x_{i}))^{\gamma}} = 0,$$

$$\hat{\alpha} = n\hat{\beta}\left[(\hat{\beta}+1)\sum_{i=1}^{n}\frac{(1-F(x_i))^{\hat{\gamma}}}{1+(\alpha-1)(1-F(x_i))^{\hat{\gamma}}}\right]^{-1}. \tag{14}$$

$$\frac{\partial[\ln L(\underline{x}|\alpha,\beta,\gamma)]}{\partial \beta} = \frac{n}{\beta} + n\ln(\alpha) + \gamma\sum_{i=1}^{n}\ln(1-F(x_i)) - \sum_{i=1}^{n}\ln\left(1+(\alpha-1)(1-F(x_i))^{\gamma}\right) = 0,$$

$$\hat{\beta} = n\left[\sum_{i=1}^{n}\ln\left(1-(1-\hat{\alpha})(1-F(x_i))^{\hat{\gamma}}\right) - \hat{\gamma}\sum_{i=1}^{n}\ln(1-F(x_i)) - n\ln\hat{\alpha}\right]^{-1}. \tag{15}$$

$$\frac{\partial[\ln L(\underline{x}|\alpha,\beta,\gamma)]}{\partial \gamma} = \frac{n}{\gamma} + \beta\sum_{i=1}^{n}\ln(1-F(x_i)) - (\beta+1)\sum_{i=1}^{n}\frac{(\alpha-1)(1-F(x_i))^{\hat{\gamma}}\ln(1-F(x_i))}{1+(\alpha-1)(1-F(x_i))^{\hat{\gamma}}} = 0,$$

We compute the component of Fisher Information matrix as

$$I_{\alpha\alpha} = \frac{\partial^2[\ln L(\underline{x}|\alpha,\beta,\gamma)]}{\partial \alpha^2} = -\frac{n\beta}{\alpha^2} + (\beta+1)\sum_{i=1}^{n}\frac{(1-F(x_i))^{2\gamma}}{\left[1-(1-\alpha)(1-F(x_i))^{\gamma}\right]^2}, \tag{16}$$

$$I_{\beta\beta} = \frac{\partial^2[\ln L(\underline{x}|\alpha,\beta,\gamma)]}{\partial \beta^2} = -\frac{n}{\beta^2}. \tag{17}$$

$$I_{\gamma\gamma} = \frac{\partial^2[\ln L(\underline{x}|\alpha,\beta,\gamma)]}{\partial \gamma^2} = -\frac{n}{\gamma^2} - \frac{(\beta+1)\sum_{i=1}^{n}(\alpha-1)(1-F(x))^{2\gamma}(\ln(1-F(x)))^2\left[(1+(1-F(x))^{\gamma}(\alpha-1)-(\alpha-1)\right]}{\left(1+((1-F(x))^{\gamma}(\alpha-1))\right)}$$

$$I_{\alpha\beta} = I_{\beta\alpha} = \frac{\partial^2[\ln L(\underline{x}|\alpha,\beta,\gamma)]}{\partial \alpha \partial \beta} = \frac{n}{\alpha} - \sum_{i=1}^{n}\frac{(1-F(x_i))^{\gamma}}{1+(\alpha-1)(1-F(x_i))^{\gamma}}. \tag{18}$$

$$I_{\alpha\gamma} = I_{\gamma\alpha} = \frac{\partial^2[\ln L(\underline{x}|\alpha,\beta,\gamma)]}{\partial \alpha \partial \gamma} = \frac{(\beta+1)\sum_{i=1}^{n}(1-F(x))^{\gamma}(\ln(1-F(x)))\left[(\alpha-1)-\left(1+(1-F(x))^{\gamma}(\alpha-1)\right)\right]}{\left(1+(1-F(x))^{\gamma}(\alpha-1)\right)}.$$

$$I_{\beta\gamma} = I_{\gamma\beta} = \frac{\partial^2[\ln L(\underline{x}|\alpha,\beta,\gamma)]}{\partial \beta \partial \gamma} = \sum_{i=1}^{n}\ln(1-F(x)) - \frac{\sum_{i=1}^{n}(\alpha-1)(1-F(x))^{\gamma}\ln(1-F(x))}{1+(1-F(x))^{\gamma}(\alpha-1)}.$$

The observed Fisher Information matrix is:

$$I(\underline{\theta}) = \begin{bmatrix} I_{\alpha\alpha} & I_{\alpha\beta} & I_{\alpha\gamma} \\ I_{\beta\alpha} & I_{\beta\beta} & I_{\beta\gamma} \\ I_{\gamma\alpha} & I_{\gamma\beta} & I_{\gamma\gamma} \end{bmatrix}.$$

The expected Fisher Information matrix is defined as $-E[I(\underline{\theta})]$.

The interval estimation and hypothesis testing for the parameters in $\underline{\theta}$, we require $(p+3)\times(p+3)$ matrix where $p$ denotes the number of parameters in base line distribution. Total observation matrix $I(\underline{\theta})$, where the elements of this matrix are calculated numerically. The estimated multivariate normal distribution $N_{p+3}\left(\theta, I(\hat{\theta})^{-1}\right)$ can be used to construct approximate confidence regions for the parameters in $\hat{\theta}$. An asymptotic confidence interval (ACI) with significance level $\gamma$ for each parameter $\theta_r$ is given by

$$ACI\left(\theta_r, 100(1-\gamma)\%\right) = \left(\hat{\theta}_r - z_{\gamma/2}\sqrt{\hat{\kappa}^{\theta_r,\theta_r}}, \hat{\theta}_r + z_{\gamma/2}\sqrt{\hat{\kappa}^{\theta_r,\theta_r}}\right), \tag{19}$$

where $\hat{\kappa}^{\theta_r,\theta_r}$ is the rth diagonal element of $I(\theta)^{-1}$ estimated at $\hat{\theta}$ and $z_{\gamma/2}$ is the quantile $1-\gamma/2$ of the standard normal distribution.

To test sub-models of GEMO family we compute the maximum likelihood values of the unrestricted and restricted log-likelihoods to construct likelihood ratio (LR) statistic. In any case, considering the partition $\theta = \left(\theta_1^T, \theta_2^T\right)^T$, tests of hypotheses of the type $H_0 : \theta_1 = \theta_1^{(0)}$ versus $H_A : \theta_1 \neq \theta_1^{(0)}$ can be performed using the LR statistic $LRT = 2\{l(\hat{\theta}) - l(\tilde{\theta})\}$, where $\hat{\theta}$ and $\tilde{\theta}$ are the estimates of $\theta$ under $H_A$ and $H_0$, respectively. Under the null hypothesis $H_0$, $LRT \xrightarrow{d} \chi_q^2$, where $q$ is the dimension of the vector $\theta_1$ of interest. The LR test rejects $H_0$ if $t > \xi_\gamma$, where $\xi_\gamma$ denotes the upper $100\gamma\%$ point of the $\chi_q^2$ distribution.

## 5. Applications

For illustration purpose, we use Weibull as baseline distribution (GEMO-W) and compare it with five-parameter Dias et al., (2016) Weibull (EMO-G-W), Weibull, Gamma and Exponential distributions. The maximization is carried out by using Broyden–Fletcher–Goldfarb–Shanno (BFGS) method (also known as a variable metric algorithm) with analytical derivatives, using R-language, R (2013).

### 5.1 Example 1

The bladder cancer is a disease in which abnormal cells multiply without control in the bladder. One of the most common type of bladder cancer recapitulates the normal histology of the urothelium and is known as transitional cell carcinoma. We take the data set used by Lee and Wang (2003) and consider the remission time (in months) of 128 bladder cancer patients. We assume that the data sets are uncensored. The total time on test (TTT) plot, Aarset (1987), is used to find whether the hazard rate is increasing or decreasing.

The TTT plot and model fit for bladder cancer data is shown in Figure 11 (a and b) respectively. The TTT curve; provide evidence that a bath tub hazard rate is adequate and the GEMO-W model gives the better fit as compared to the other fitted models for bladder cancer patients' data.

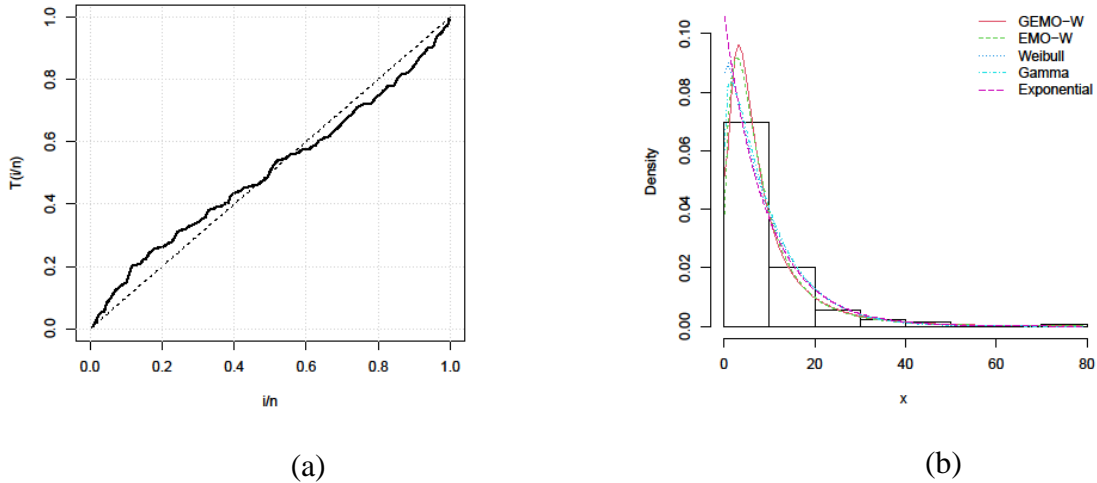

(a)            (b)

**Figure 11:** (a) The TTT plot (b) Histogram along with the fitted pdfs, for bladder cancer data

In Table 1, the MRL and MPL are calculated for the cancer patients by using GEMO-W. The MRL shows the expected value of the remaining lifetime after different percentiles time points. As the remission time increases for the cancer patients the mean residual life also increases. The MPL corresponds to the mean time elapsed since the failure of $T$, given that $T \leq t$. In real life situations usually one needs this function to access the lifetime of the system. As the remission time increases for the cancer patients the mean past lifetime also increases.

**Table 1**: Mean Residual Life and Mean Past Lifetime for Cancer Patients

| Percentiles | 0.1 | 0.25 | 0.5 | 0.75 | 0.9 | 0.95 | 0.99 |
|---|---|---|---|---|---|---|---|
| Time (Month) | 1.67 | 3.36 | 6.65 | 11.98 | 19.94 | 26.26 | 45.47 |
| MRL | 8.67 | 8.90 | 9.16 | 10.69 | 13.16 | 15.07 | 20.28 |
| MPL | 0.73 | 1.50 | 3.27 | 6.85 | 13.24 | 18.82 | 36.95 |

## 5.2 Example-2

In this example, we use the data of the strength of glass fibers of length 1.5 cm taken from the National Physical Laboratory in England Smith and Naylor (1987). The TTT plot and model fit for the strength of glass fiber data is shown in Figure 12 (a and b) respectively. The TTT curve is concave, which shows that a monotonic increasing hazard rate is adequate and the GEMO-W model gives a better fit as compared to the other fitted models for this data set.

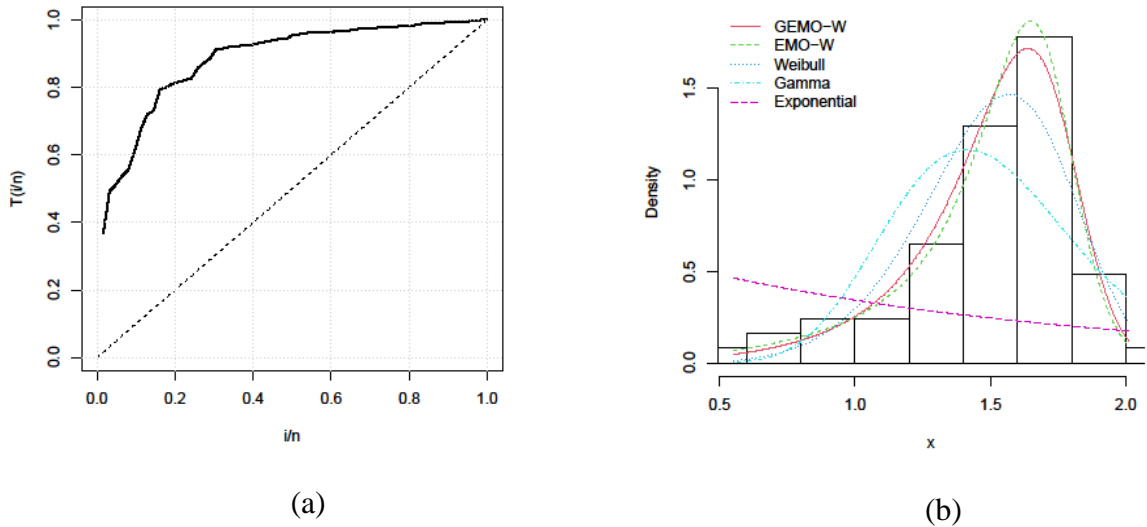

(a)                                                (b)

**Figure 12:** (a) The TTT plot (b) Histogram along with the fitted pdfs for strength of glass fibers data

In Table 2 the MRL and MPL are calculated for the strength of glass fiber by using GEMO-W. As strength increases for the glass fiber the mean residual life decrease and mean past lifetime increases.

**Table 2**: Mean Residual Life and Mean Past Lifetime for Glass Fiber

| Percentiles | 0.1  | 0.25 | 0.5  | 0.75 | 0.9  | 0.95 | 0.99 |
|-------------|------|------|------|------|------|------|------|
| Strength    | 1.05 | 1.39 | 1.60 | 1.69 | 1.81 | 1.85 | 2.00 |
| MRL         | 0.52 | 0.26 | 0.14 | 0.11 | 0.08 | 0.07 | 0.05 |
| MPL         | 0.26 | 0.29 | 0.29 | 0.30 | 0.36 | 0.38 | 0.51 |

The results of parameter estimation and model selection criteria have been reported in Table 3. The standard errors are given in the parentheses below the estimates. For model selection, we report the following criteria: the Negative log Likelihood (LL), Akaike Information Criterion (AIC), Kolmogorov-Smirnov (KS) statistic and Anderson-Darling (AD) Statistic. The model with the highest negative log-likelihood value and the smallest AIC, KS and AD values indicates the better model. The standard errors for $\beta$ and $\lambda$ are smaller than their ML estimates which show greater stability in their estimation. Although the standard errors of $\alpha, \theta$ and $\gamma$ are larger than their ML estimates showing more variability in their estimation but from Table 3 the values of Negative-LL, AIC, KS statistic and AD statistic for the GEMO-W are comparatively better than the other models in both examples. Hence, our proposed Weibull model (GEMO-W) can be chosen as the best model for both data sets.

## 6. Conclusions

We proposed a generalized form of Marshal-Olkin family of distributions named as GEMO. This distribution included three additional parameters which aided to model the reliability/lifetime data with high or low kurtosis, heavy tailed and right sided skewness. The GEMO distribution incorporated with Weibull as a baseline distribution (GEMO-W) performed better than EMO-G-W, Weibull, Gamma and Exponential distributions for modeling data sets in the afore mentioned situations. The comparative performance of GEMO-W was observed to give an adequate "good fit" on the real life data sets based on cancer patients and glass fiber. It is a valuable new contribution to the pool of current development of Marshall-Olkin models.

**Acknowledgments:**

The authors would like to thank the two referees and the editor for their comments which greatly improved this paper. The authors would also like to thank Prof. B. C. Arnold, for his valuable suggestion to expand this research article.

**Table 3**: Estimation of the parameters, the corresponding SEs and the statistics: Log Likelihood, AIC, Kolmogorov-Smirnov (KS) and Anderson-Darling (AD) criteria

| Example-1 Cancer Patients | $\hat{\alpha}$ | $\hat{\beta}$ | $\hat{\lambda}$ | $\hat{\theta}$ | $\hat{\gamma}$ | Negative LL | AIC | KS Statistic | AD Statistic |
|---|---|---|---|---|---|---|---|---|---|
| GEMO-W | 25.5629 (43.4078) | 0.2846 (0.1915) | 0.5946 (0.2368) | 3.6174 (62.2268) | 4.0532 (41.4933) | -409.3703 | 828.7479 | 0.0306 | 0.0951 |
| EMO-G-W | 3.0110 (30.4162) | 1.9956 (0.9577) | 0.3492 (0.1493) | 2.6080 (75.2204) | -23.6085 (67.0649) | -409.4687 | 828.9369 | 0.0309 | 0.0974 |
| Weibull | 1 --- | 1 --- | 1.0477 (0.0676) | 9.5600 (0.8529) | 1 --- | -414.0869 | 832.1738 | 0.0700 | 0.9578 |
| Gamma | 1 --- | 1 --- | 1.1725 (0.1308) | 0.1252 (0.0173) | 1 --- | -413.3678 | 830.7356 | 0.0732 | 0.7706 |
| Exponential | 1 --- | 1 --- | 0.1068 (0.0094) | 1 | 1 --- | -414.3419 | 830.6838 | 0.0846 | 1.1736 |
| **Example-2 Glass Fiber** | | | | | | | | | |
| GEMO-W | 31.2329 (57.4598) | 2.5066 (2.4793) | 3.1095 (1.3128) | 2.1741 (23.0724) | 6.7698 (22.4060) | -6.2412 | 20.4824 | 0.0829 | 0.2705 |
| EMO-G-W | 17.0719 (217.178) | 0.6260 (0.1608) | 4.7370 (1.5491) | 2.3774 (6.3890) | -18.1544 (26.5353) | -6.2600 | 22.5201 | 0.0864 | 0.2859 |
| Weibull | 1 --- | 1 --- | 6.3269 (0.6665) | 1.6110 (0.0337) | 1 --- | -11.0804 | 26.1608 | 0.1524 | 1.3257 |
| Gamma | 1 --- | 1 --- | 18.0670 (3.2154) | 12.0849 (2.1809) | 1 --- | -22.0264 | 48.0527 | 0.2239 | 3.3871 |
| Exponential | 1 --- | 1 --- | 0.6689 (0.0849) | 1 --- | 1 --- | -86.9318 | 175.8636 | 0.4185 | 18.3791 |